\magnification=1250
\vsize = 25truecm
\hsize = 18truecm
\hoffset = -1truecm
\voffset = -.5truecm

\font\tengoth=eufm10 at 10 pt \font\sevengoth=eufm10 at 7 pt
\font\fivegoth=eufm10 at 5 pt \newfam\gothfam
\textfont\gothfam=\tengoth \scriptfont\gothfam=\sevengoth
\scriptscriptfont\gothfam=\fivegoth 
\newcount\subsecno
\newcount\secno
\newcount\prmno
\newif\ifnotfound
\newif\iffound

\def\namedef#1{\expandafter\def\csname #1\endcsname}
\def\nameuse#1{\csname #1\endcsname}

\long\def\ifundefined#1#2#3{\expandafter\ifx\csname
  #1\endcsname\relax#2\else#3\fi}
\def\hwrite#1#2{{\let\the=0\edef\next{\write#1{#2}}\next}}

\toksdef\ta=0 \toksdef\tb=2
\long\def\leftappenditem#1\to#2{\ta={\\{#1}}\tb=\expandafter{#2}%
                                \edef#2{\the\ta\the\tb}}
\long\def\rightappenditem#1\to#2{\ta={\\{#1}}\tb=\expandafter{#2}%
                                \edef#2{\the\tb\the\ta}}

\def\lop#1\to#2{\expandafter\lopoff#1\lopoff#1#2}
\long\def\lopoff\\#1#2\lopoff#3#4{\def#4{#1}\def#3{#2}}

\def\ismember#1\of#2{\foundfalse{\let\given=#1%
    \def\\##1{\def\next{##1}%
    \ifx\next\given{\global\foundtrue}\fi}#2}}

\def\section#1{\medbreak
               \global\def\currenvir{section}
               \global\advance\secno by1\global\prmno=0\global\subsecno=0
               {\bf \ind\number\secno. {#1}$.-$}
               }

\def\subsection{\global\def\currenvir{subsection}
                \global\advance\prmno by1\global\subsecno=0
                \ind{\bf \number\secno.\number\prmno. }}
\def\rem{\global\def\currenvir{rem}
                \global\advance\prmno by1\global\subsecno=0
                \medskip\ind{\bf Remark \number\secno.\number\prmno$.-$ }}
\def\ex{\global\def\currenvir{ex}
                \global\advance\prmno by1\global\subsecno=0
                \medskip\ind{\bf Example \number\secno.\number\prmno$.-$ }}
\def\formule{\global\def\currenvir{formule}
\ifnum\prmno=0\global\advance\prmno by1
{\number\secno.\number\prmno}\else
\global\advance\subsecno by1
 {\number\secno.\number\prmno.\number\subsecno}\fi}

\def\proclaim#1{\global\advance\prmno by 1\global\subsecno=0
                {\bf #1 \the\secno.\the\prmno$.-$ }}

\long\def\th#1 \enonce#2\endth{%
   \medbreak\proclaim{#1}{\it #2}\global\def\currenvir{th}\medskip}


\def\isinlabellist#1\of#2{\notfoundtrue%
   {\def\given{#1}%
    \def\\##1{\def\next{##1}%
    \lop\next\to\za\lop\next\to\zb%
    \ifx\za\given{\zb\global\notfoundfalse}\fi}#2}%
    \ifnotfound{\immediate\write16%
                 {Warning - [Page \the\pageno] {#1} No reference found}}%
                \fi}%
\def\ref#1{\ifx\labellist\empty{\immediate\write16
                 {Warning - No references found at all.}}
               \else{\isinlabellist{#1}\of\labellist}\fi}

\def\newlabel#1#2{\rightappenditem{\\{#1}\\{#2}}\to\labellist}
\def\labellist{}

\def\label#1{%
  \def\given{th}%
  \ifx\given\currenvir%
    {\hwrite\lbl{\string\newlabel{#1}{\number\secno.\number\prmno}}}\fi%
  \def\given{section}%
  \ifx\given\currenvir%
    {\hwrite\lbl{\string\newlabel{#1}{\number\secno}}}\fi%
   \def\given{formule}%
  \ifx\given\currenvir%
{
\ifnum\prmno=1\hwrite\lbl
{
\string\newlabel{#1}{\number\secno.\number\prmno}
}
\else\hwrite\lbl
{
\string\newlabel{#1}{\number\secno.\number\prmno.\number\subsecno}
}\fi}
\fi%
  \def\given{subsection}%
  \ifx\given\currenvir%
    {\hwrite\lbl{\string\newlabel{#1}{\number\secno.\number\prmno}}}\fi%
\def\given{rem}%
  \ifx\given\currenvir%
    {\hwrite\lbl{\string\newlabel{#1}{\number\secno.\number\prmno}}}\fi%
\def\given{ex}%
  \ifx\given\currenvir%
    {\hwrite\lbl{\string\newlabel{#1}{\number\secno.\number\prmno}}}\fi%
  \def\given{subsubsection}%
  \ifx\given\currenvir%
  {\hwrite\lbl{\string%
    \newlabel{#1}{\number\secno.\number\subsecno.\number\subsubsecno}}}\fi
  \ignorespaces}

\newwrite\lbl

\def\openall{\openout\lbl=\jobname.lbl}
\def\closeall{\closeout\lbl}

\newread\testfile
\def\lookatfile#1{\openin\testfile=\jobname.#1
    \ifeof\testfile{\immediate\openout\nameuse{#1}\jobname.#1
                    \write\nameuse{#1}{}
                    \immediate\closeout\nameuse{#1}}\fi%
    \immediate\closein\testfile}%

\def\begin{\lookatfile{lbl}
           \input\jobname.lbl
           \openall}
\let\bye\end
\def\end{\closeall\bye}


\mathcode`A="7041 \mathcode`B="7042 \mathcode`C="7043
\mathcode`D="7044 \mathcode`E="7045 \mathcode`F="7046
\mathcode`G="7047 \mathcode`H="7048 \mathcode`I="7049
\mathcode`J="704A \mathcode`K="704B \mathcode`L="704C
\mathcode`M="704D \mathcode`N="704E \mathcode`O="704F
\mathcode`P="7050 \mathcode`Q="7051 \mathcode`R="7052
\mathcode`S="7053 \mathcode`T="7054 \mathcode`U="7055
\mathcode`V="7056 \mathcode`W="7057 \mathcode`X="7058
\mathcode`Y="7059 \mathcode`Z="705A
\def\spacedmath#1{\def\packedmath##1${\bgroup\mathsurround =0pt##1\egroup$}
\mathsurround#1
\everymath={\packedmath}\everydisplay={\mathsurround=0pt}}
\def\nospacedmath{\mathsurround=0pt
\everymath={}\everydisplay={} } \spacedmath{2pt}
\def\qfl#1{\normalbaselines{\baselineskip=0pt
\lineskip=10truept\lineskiplimit=1truept}\nospacedmath\smash {\mathop{\hbox to
6truemm{\rightarrowfill}}
\limits^{\scriptstyle#1}}}

\def\phfl#1#2{\normalbaselines{\baselineskip=0pt
\lineskip=10truept\lineskiplimit=1truept}\nospacedmath\smash {\mathop{\hbox to
8truemm{\rightarrowfill}}
\limits^{\scriptstyle#1}_{\scriptstyle#2}}}
\def\diagram#1{\def\normalbaselines{\baselineskip=0truept
\lineskip=10truept\lineskiplimit=1truept}   \matrix{#1}}

\def\Pic{\mathop{\rm Pic}\nolimits}
\def\pc#1{\tenrm#1\sevenrm}
\def\up#1{\raise 1ex\hbox{\smallf@nt#1}}
\def\tx{\kern-1.5pt -}
\def\cqfd{\kern 2truemm\unskip\penalty 500\vrule height 4pt depth 0pt width
4pt\medbreak} 
\def\virg{\raise .4ex\hbox{,}}
\def\decale#1{\smallbreak\hskip 28pt\llap{#1}\kern 5pt}
\def\no{n\up{o}\kern 2pt}
\def\ind{\par\hskip 0cm\relax}

\parindent=0cm

\def\rond{\kern 1pt{\scriptstyle\circ}\kern 1pt}

\baselineskip15pt
\overfullrule=0pt

\def\boxit#1#2{
\setbox1=\hbox{\kern#1{#2}\kern#1}
\dimen1=\ht1 \advance\dimen1 by #1 \dimen2=\dp1 \advance\dimen2 by #1
\setbox1=\hbox{\vrule height\dimen1 depth\dimen2\box1\vrule}
\setbox1=\vbox{\hrule\box1\hrule}
\advance\dimen1 by .4pt \ht1\dimen1
\advance\dimen2 by .4pt \dp1\dimen2 \box1\relax}

\def\varalpha{\rlap{\raise 1pt\hbox{$\scriptscriptstyle /$}}\alpha}

\def\SL{\mathop{\bf SL}\nolimits}

\let\ra\rightarrow

\def\Id{\hbox{\rm Id}}
\def\isom{\buildrel\sim\over\ra}

\def\Z{\Bbb Z}\def\C{\Bbb C}\def\F{\Bbb F}
\def\spec{\mathop{\rm spec}\nolimits}

\def\fhd#1#2{\nospacedmath\smash{\mathop{\hbox to 8mm{\rightarrowfill}}
\limits^{\scriptstyle#1}_{\scriptstyle#2}}}
\def\fhg#1#2{\nospacedmath\smash{\mathop{\hbox to 8mm{\leftarrowfill}}
\limits^{\scriptstyle#1}_{\scriptstyle#2}}}
\def\fvb#1#2{\nospacedmath\llap{$\scriptstyle #1$}\left\downarrow
\vbox to 6mm{}\right.\rlap{$\scriptstyle#2$}}
\def\fvh#1#2{\nospacedmath\llap{$\scriptstyle #1$}\left\uparrow
\vbox to 6mm{}\right.\rlap{$\scriptstyle#2$}}
 
\def\FFhd#1#2{\nospacedmath\smash{\mathop{\hbox to 8mm{\hfil$\Longrightarrow$\hfil}}
\limits^{\scriptstyle#1}_{\scriptstyle#2}}}
\def\FFhg#1#2{\nospacedmath\smash{\mathop{\hbox to 8mm{$\Longleftarrow$}}
\limits^{\scriptstyle#1}_{\scriptstyle#2}}}
\def\FFvb#1#2{\nospacedmath\llap{$\scriptstyle #1$}\left\Downarrow
\vbox to 6mm{}\right.\rlap{$\scriptstyle#2$}}

\def\Fhd{\fhd{}{}}

\def\epi{\rightarrow \kern -3mm\rightarrow }
\def\rond{\kern 1pt{\scriptstyle\circ}\kern 1pt}
\def\ra{{\rightarrow}}

\def\\{{\setminus}}

\def\Hom{\mathop{\rm Hom}\nolimits}
\def\hom{\mathop{{\cal H}om}\nolimits}

\def\P{{\bf P}}
\def\O{{\cal O}}

\let\ra\rightarrow

\def\Bbb{{\bf}}
\def\spec{{\rm Spec}}
\begin

\null
\centerline{{\bf Linearization of group stack actions and the Picard group}} 

\centerline{{\bf of the moduli
of $\SL_r/\mu_s$-bundles on a curve}} 
\medskip

\centerline{Yves {\pc LASZLO}
\footnote{\parindent 0.5cm($\dagger$)}{\sevenrm Partially
supported by the European HCM Project ``Algebraic Geometry in Europe"
(AGE).}} 

\def\a{{\cal A}}\def\b{{\cal B}}\def\c{{\cal C}}\def\tp{{\tilde p}}\def\X{{\cal X}}        
\def\D{{\cal X}}\def\C{{\cal G}}\def\L{{\cal L}}\def\E{{\cal Y}}\def\m{{m_\C}}\def\M{{\cal M}}
\def\J{{\cal J}}\def\p{{\pi}}\def\k{{\bf k}}\def\P{{\cal P}}\def\T{{\cal T}}\def\H{{\cal H}}
\def\U{{\cal U}}\def\d{{\cal D}}\def\F{{\cal F}}
\let\scp=\scriptstyle
{\bf Introduction} 
\medskip
Let $G$ be a complex semi-simple group and $\tilde G\epi G$ the universal covering. Let $\M_G$
(resp. $\M_{\tilde G}$) be the moduli stack of $G$-bundles over a curve $X$ of degree
$1\in\pi_1(G)$ (resp. of ${\tilde G}$-bundles. In [B-L-S], we have studied the link between the
groups $\Pic(\M_G)$ and $\Pic(\M_{\tilde G}$), the later being well understood thanks to [L-S].
In particular, it has been possible to give a complete description in the case where $G={\bf
PSL}_r$ but not in the case $\SL_r/\mu_s,\ s\mid r$, although we were able to give partial
results. The reason was that we did not have at our disposal the technical background to study
the morphism $\M_{\tilde G}\ra\M_G$. It turns to be out that it is a torsor under some group
stack, not far from a Galois \'etale cover in the usual schematic picture. Now, the descent
theory of Grothendieck has been adapted to the set-up of fpqc morphisms of stacks  in [L-M] and
gives the theorem
\ref{theo} in the particular case of a morphism which
is torsor under a group stack. We then used this technical result to determine the exact structure
of
$\Pic(\M_G)$ where
$G=\SL_r/\mu_s$ (theorem
\ref{theo-pic}).
\medskip I would like to thank L. Breen to have taught me both the notion of torsor and
linearization of a vector bundle in the set-up of group-stack action and for his comments on a
preliminary version of this paper. 
\medskip 
{\bf Notation}
\medskip
Throughout this paper, all the stacks will be implicitely assumed to be algebraic over a fixed base
scheme and the morphisms locally of finite type. We fix once for all a projective, smooth,
connected genus $g$ curve $X$ over an algebraically closed field $\k$ and a closed point $x$ of
$X$. For simplicity, we assume $g>0$ (see remarks \ref{rema} and \ref{rem} for
the case of
${\bf P}^1$). The Picard stack parametrizing families of line bundles of degree
$0$  on $X$ will be denoted by $\J(X)$ and the jacobian variety of $X$ by $JX$.
If $G$ is an algebraic group over $\k$, the quotient stack $\spec(\k)/G$ (where
$G$ acts trivialy on $\spec(\k))$ whose category over a
$\k$-scheme $S$ is the category of $G$-torsors (or $G$-bundles) over $S$ will be denoted by $BG$.
If $n$ is an integer and $A=\J(X), JX$ or $BG_m$ we denote by $n_A$  the $n^{\rm\scp th}$-power
morphism $a\longmapsto a^n$. We denote by $\J_n$ (resp. $J_n$) the $0$-fiber $A\times_A\spec(\k)$
of $n_A$ when $A=\J(X)$ (resp. $A=JX$), wh

\section{Generalities}  Following [Br], for
any diagram $$ A\fhd{h}{}B\matrix{\fhd{g}{}\cr\Uparrow\rlap{$\scp\lambda$}\cr
\scp f\cr\Fhd}C\fhd{l}{}D$$ of $2$-categories, we'll denote by $l*\lambda:\ l\rond f\Rightarrow
l\rond g$ (resp. $\lambda*h:\ f\rond h\Rightarrow g\rond h$) the $2$-morphism deduced from
$\lambda$.
\subsection For the convenience of the reader, let us prove a simple formal lemma which will be usefull in the 
section \ref{simplicial}. Let $\a,\b,\c$ be three $2$-categories, a $2$-commutative
diagram $$\matrix{
&&\c\cr
&\llap{${}^{\delta_0}$}\nearrow&\fvh{}{d_0}\cr
\a&\fhd{f}{}&\b\cr
&\llap{${}_{\delta_1}$}\searrow&\fvb{}{d_1}\cr
&&\c\cr}\leqno{(\formule)}$$\label{strict-categorie}
and a $2$-morphism $\mu:\ \delta_0\Rightarrow \delta_1$.
\th Lemma 
\enonce Assume that $f$ is an equivalence. There exists a unique $2$-morphism $$\mu*f^{-1}:\
d_0\Rightarrow d_1$$ such that $(\mu*f^{-1})*f=\mu$.
\endth\label{lemme-categorie}
{\it Proof}: let $\epsilon_k, k=0,1$ the $2$-morphism $d_k\rond f\Rightarrow \delta_k$. Let $b$ be an object of $\b$. Pick an object $a$ of $\a$ and an isomorphism
$\alpha:\ f(a)\isom b$. Let $\varphi_\alpha :\ d_0(b)\isom d_1(b)$ be the unique isomorphism making
the diagram
$$\matrix{
\delta_0(a)&\fhd{\epsilon_0(a)}{}&d_0\rond f(a)&\fhd{d_0(\alpha )}{}&d_0(b)\cr
\fvb{\mu_a}{}&&&&\fvb{}{\varphi_\alpha }\cr
\delta_1(a)&\fhd{\epsilon_1(a)}{}&d_1\rond f(a)&\fhd{d_1(\alpha )}{}&d_1(b)\cr
}$$ commutative. We have to show that $\varphi_\alpha $ does not depend on $\alpha $ but only
on $b$. Let $\alpha ':\ f(a')\isom b$ be another isomorphism. There exists a unique isomorphism
$\iota:\ a'\isom a$ such that $\alpha\rond f(\iota)=\alpha'$. The one has the equality
$\varphi_{\alpha'}=d_1(\alpha )\rond\Phi\rond d_0(\alpha)^{-1}$ where 
$$\Phi=[d_1\rond
f(\iota)]\rond\epsilon_1(a')\rond \mu_{a'}\rond\epsilon_0(a')^{-1}\rond[d_0\rond f(\iota)]^{-1}.$$ 
The functoriality of $\epsilon_i$ and $\mu$ ensures that one has the equalities
$$d_k\rond f(\iota)\rond\epsilon_k(a')=\epsilon_k(a)\rond\delta_k(\iota)$$ and
$$\mu_{a}=\delta_1(\iota)\rond\mu_{a'}\rond\delta_0(\iota)^{-1}.$$ This shows that
$$\Phi=\epsilon_1(a)\rond\mu_a\rond\epsilon_0(a)^{-1}$$ which proves 
the equality $\varphi_\alpha =\varphi_{\alpha '}$. We can therefore define $\mu_b$ 
to be the isomorphism
$\varphi_\alpha $ for one isomorphism $\alpha:\ f(a)\isom b$. One checks that the construction
is functorial in $b$ and the lemma follows.\cqfd

\section{Linearizations of line bundles on stacks} Let us first recall following [Br] the notion
of torsor in the stack context.
\subsection Let $f:\ \D\ra\E$ be a faithfully flat morphism of stacks.  Let us assume that an
algebraic $gr$-stack $\C$ acts on $f$ (the product of $\C$ is denoted by $\m$ and the unit object
by $1$). Following [Br], this means that there exists a 1-morphism of $\E$-stacks $m:\
\C\times\D\ra\D$ and a $2$-morphism $\mu:\ m\rond (\m\times \Id_\D)\Rightarrow m\rond
(\Id_\C\times m)$  such that the obvious associativity condition (see the diagram (6.1.3) of
[Br]) is satisfied and such that there exists a $2$-morphism $\epsilon:\ m\rond(1\times
\Id_\D)\Rightarrow \Id_\D$ which is compatible to $\mu$ in the obvious sense (see (6.1.4) of
[Br]). 

\rem To say that $m$ is a morphism of $\E$-stacks means that the diagram
$$\matrix{\C\times\D&\fhd{m}{}&\D\cr\searrow&&\swarrow\cr&\E\cr}$$ is $2$-commutative. In other
words, if we denote for simplicity the image of a pair of objects $m(g,x)$ by $g.x$, this means
that there exists a functorial isomorphism $\iota_{g,x}:\ f(g.x)\ra f(x)$.\label{Y-morphism}
 
\subsection\label{equi} Suppose that $\C$ acts on another such $f':\D'\ra\E$. A morphism
$p:\ \D'\ra\D$ will be said equivariant if there exists a  $2$-morphism 
$$q:\ m\rond(\Id\times p)\Rightarrow p\rond m'$$ which is compatible to $\mu$ (as in
[Br] (6.1.6)) and $\epsilon$ (which is implicit in [Br]) in the obvious sense.

\th Definition 
\enonce With the above notations, we say that $f$ (or $\D$) is a $\C$-torsor over
$\E$ if the morphism $pr_2\times m:\ \C\times \D\ra\D\times_\E\D$ is an isomorphism (of stacks)
and  the geometrical fibers of $f$ are not empty. 
\endth

\rem In down to earth terms, this means that if $\iota:\ f(x)\ra f(x')$ is an isomorphism in
$\E$ ($x,x'$ being objects of $\D$), there exist an object $g$ of $\C$ and a unique isomorphism
$(x,g.x)\fhd{\sim}{}(x,x')$ which induces $\iota$ thanks to $\iota_{g,x}$ (cf. \ref{Y-morphism}).

\ex If $\M_X(G_m)$ is the
Picard stack of $X$, the morphism $\M_X(G_m)\ra\M_X(G_m)$ of multiplication by $n\in\Z$ is a
torsor under $B\mu_n\times J_n(X)$ (cf. (\ref{torseur})).

\subsection Let a $\L$ be a line bundle  on $\D$. By definition, the data $\L$
 is equivalent to te data of a morphism $l:\ \D\ra BG_m$ (see [L-M],prop. 6.15). If $\L,\L'$ are
$2$ line bundles on $\X$ defined by $l,l'$, we will view an isomorphism 
$\L\isom\L'$ as a $2$-morphism $l\Rightarrow l'$.\label{convention-fibre}
 
\th Definition 
\enonce A $\C$-linearization of $\L$ is a $2$-morphism $\lambda:\ l\rond
m\Rightarrow l\rond pr_2$ such that the two diagrams of $2$-morphisms $$\matrix{ l\rond m\rond
(\m\times \Id_\D)&\buildrel{l*\mu}\over{\Longrightarrow}&l\rond m\rond (\Id_\C\times m)\cr
\hfill\Big\Downarrow{\scriptstyle\lambda*(\m\times
\Id_\D)}&&\hfill\Big\Downarrow{\scriptstyle\lambda*(\Id_\C\times m)}\cr l\rond pr_2\rond(\m\times
\Id_\D)=l\rond pr_2\rond pr_{23}&\buildrel \lambda*pr_{23}\over{\Longleftarrow}&l\rond
pr_2\rond(\Id_\C\times m)=l\rond m\rond pr_{23}\cr }\leqno(\formule)$$\label{lineA} and $$
\matrix{ l\rond m\rond(1\times \Id_\D)&\buildrel l*\epsilon\over{\Longrightarrow}&l\cr
\FFvb{\lambda*(1\times \Id_\D)}{}&&\parallel\cr l&=&l\cr
}\leqno(\formule)$$\label{lineB}(stritly) commutes. 
\endth

\rem\label{rem-line} In $g_1,g_2$ are objects of $\C$ and $d$ of $\D$, the commutativity of the
diagram (\ref{lineA}) means that the diagram 
$$\matrix{ 
\L_{(g_1.g_2)x}&\fhd{\sim}{}&\L_{g_1(g_2.x)}\cr
\fvb{}{\wr}&&\fvb{}{\wr}\cr 
\L_x&\fhg{\sim}{}&\L_{g_2.x}\cr}$$ is commutative and the commutativity of (\ref{lineB}) that the 
two
isomorphisms $\L_{1.x}\simeq \L_x$ defined by the linearization $\lambda$ and $\epsilon$
respectively are the same.

\section{An example} Let me recall that a closed point $x$ of $X$ has been fixed. Let $S$ be a
$\k$-scheme. The $S$-points of the jacobian variety of $X$ are by definition isomorphism classes of
line bundles on $X_S$ together whith a trivialization along $\{x\}\times S$ (such a pair will be
called a rigidified line bundle). For the covenience of the reader, let me state this well known
lemma which can be founf in SGA4, exp. XVIII, (1.5.4) 
\th Lemma
\enonce The Picard stack $\J(X)$ is canonically isomorphic (as a $\k$-group stack) to $JX\times
BG_m$.
\endth\label{torseur}
{\it Proof}: let $f:\ \J(X)\ra JX\times BG_m$ be the morphim which associates

-to the line bundle $L$ on $X_S$ the pair 
$L\otimes L_{\mid \{x\}\times S}^{-1},L_{\mid \{x\}\times S}$ (thought as an object of $JX\times
BG_m$ over $S$);

-to an isomorphism $L\isom L'$ on $X_S$ its restriction to $\{x\}\times S$.

Let $f':\ JX\times BG_m\ra\J(X)$ be the morphism which associates

-to the pair $(L,V)$ where $L$ is a rigidified bundle on $X_S$ and $V$ a line bundle on $S$
(thought as an object of $JX\times BG_m$ over $S$), the line bundle $L\otimes_{X_R}V$;

-to an isomorphism $(l,v):\ (L,V)\isom (L',V')$ the tensor product $l\otimes_{X_S} v$.

The morphisms $f$ and $f'$ are (quasi)-inverse each other and are morphisms of $\k$-stacks.\cqfd

We will identify from now $\J(X)$ and $JX\times BG_m$. Let $\L$ (resp. $\P$ and $\T$) be the
universal bundle on $X\times\J(X)$ (resp. on $X\times JX$ and $BG_m$) and let $\Theta=(\det
R\Gamma\P)^{-1}$ be the theta line bundle on $JX$. The isomorphism $\L\isom\P\otimes\T$ yields
an isomorphism $$\det
R\Gamma\L^n(m.x)\isom\Theta^{-n^2}\otimes\T^{(m+1-g)}.\leqno({\formule})$$\label{iso-det}

 \section{Descent of $\C$-line bundles} \label{simplicial} The object of this section is to prove
the following statement

\th Theorem
\enonce Let $f:\D\ra\E$ a $\C$-torsor as above. Let $\Pic^\C(\D)$ be the group of isomorphism
classes of $\C$-linearized line bundles on $\D$. 

Then, the pull-back morphism $f^*:\
\Pic(\E)\isom\Pic^\C(\D)$ is an isomorphism.
\endth\label{theo}
 The descent theory of Grothendieck has been adapted in the case of algebraic $1$-stacks in
[L-M], essentially in the proposition (6.23).Let $\D_\bullet\ra\E$ be the (augmented) simplicial
complex of stacks coskeleton of $f$ (as defined in [De] (5.1.4) for instance). By proposition
(6.23) of [L-M], one just has to construct a cartesian $\O_{D_\bullet}$-module $\L_\bullet$ such
that $L_0$ is the $\O_\D$-module $\L$ to prove the theorem. The $n$-th piece $\D_n$ is inductively
defined  by $\D_0=\D,\ \D_n=\D\times_\E\D_{n-1}$ for $n>0$. Let $p_n:\ \D_n\ra\D$ be the
projection on the first factor. It is the simplicial morphism associated to the map $$\tp_n:\
\left\{\matrix{\Delta_0&\ra&\Delta_n\cr 0&\longmapsto&0\cr}\right.$$ We define $\L_n$ by the
morphism $$l_n:\ \D_n\fhd{p_n}{}\D\fhd{l}{}BG_m\leqno{(\formule)}.$$ 
\subsection\label{rel-iso} Let $\delta_i$ (resp. $s_j$) be the
face (resp. degeneracy) operators (see [De] (5.1.1) for instance) (by abuse of notation, we use the
same notation for $\delta_j,s_j$ and their image by $\D_\bullet)$. 
The category $(\Delta)$ is generated by the
face and degeneracy operators with the following relations (see for instance the proposition
VII.5.2 page 174 of [McL]) $$\matrix{\delta_i\rond\delta_j&=&\delta_{j+1}\rond\delta_i&\kern
.5cm&i\leq j\cr}\leqno{(\formule)}$$\label{A} 
$$\matrix{s_j\rond s_i&=&s_i\rond s_{j+1}&\kern
.5cm&i\leq j\cr}\leqno{(\formule)}$$\label{B} 
$$\left\{\matrix{ s_j\rond\delta_i&=&\delta_i\rond
s_{j-1}&\kern .5cm&i<j\cr &=&1&\kern .5cm&i=j,i=j+1\cr
&=&\delta_{i-1}\rond s_j&\kern .5cm&i>j+1.\cr
}
\right.\leqno{(\formule)}$$\label{C}
Therefore, the data of a cartesian $\O_{\D_\bullet}$-module $\L_\bullet$ is equivalent to the
data of isomorphisms $\alpha _j:\ \delta_j^*\L_n\isom\L_{n+1},\ j=0,\ldots,n+1$ and
$\beta _j:\ s_j^*\L_{n+1}\isom\L_{n},\ j=0,\ldots,n$ (where $n$ is a non negative integer) which
are compatible with the relations \ref{A}, \ref{B} and \ref{C}. 

Let $n$ be a non negative integer. 
\subsection We have first to define for $j=0,\ldots,n+1$ an isomorphism $\alpha _j:\
\delta_j^*\L_n\isom\L_{n+1}$. The line bundle $\delta_j^*\L_n$ is defined by the morphism $l\rond
p_n\rond\delta_j:\ \D_{n+1}\ra BG_m$ and $\tp_n\rond\delta_j$ is associated to the map
$$\left\{\matrix{\Delta_0&\ra&\Delta_{n+1}\cr 0&\longmapsto&\delta_j(0)\cr}\right.$$ If $j\not=0$,
one has therefore $\tp_n\rond\delta_j=\tp_{n+1}$ and $\delta_j^*\L_n=L_{n+1}$. We define $\alpha
_j$ by the identity in this case. 

Suppose now that $j=0$. Let $\pi_n:\ \D_n\ra\D_1$ be the projection on the $2$ first factors
(associated to the canonical inclusion $\Delta_1\hookrightarrow\Delta_n$. The commutativity of the
$2$ diagrams
$$\matrix{
\D_{n+1}&\fhd{\delta_0}{}&\D_{n}\cr\fvb{}{\pi_{n+1}}&&\fvb{p_{n}}{}\cr\D_1&\fhd{\delta_0}{}&\D\cr}
\ {\rm and}\
\matrix{
\D_{n+1}&\fhd{p_{n+1}}{}&\D\cr\fvb{}{\pi_{n+1}}&&\fvh{\delta_1}{}{}\cr\D_1&=&\D_1\cr}$$ allows to
reduce the problem to the construction of an isomorphism $$\delta_0^*\L\isom \delta_1^*\L\ {\rm where}\
\delta_i,\D_1\ra\D\ i=0,1$$ are the face morphisms or, what is amounts to the same, to the construction
of a $2$-morphism $\nu:\ l\rond \delta_0\Rightarrow l\rond \delta_1$ (the morphism $\alpha _j$
will be $\alpha _j=\nu*\pi_{n+1}$). Now the diagram $$\matrix{ &&BG_m\cr &\llap{${}^{l\rond
m}$}\nearrow&\fvh{}{l\rond \delta_0}\cr \C\times\X&\fhd{pr_2\times m}{}&\X\times_\E\X\cr
&\llap{${}_{l\rond pr_2}$}\searrow&\fvb{}{l\rond \delta_1}\cr
&&BG_m\cr}\leqno{(\formule)}$$ is strict commutative and $pr_2\times m$ is an equivalence by the
definition of a torsor. By the lemma \ref{lemme-categorie}, the $2$-morphism $\lambda$ induces a
canonical $2$-morphism $\lambda*(pr_2\times m)^{-1}:\ l\rond \delta_0\Rightarrow l\rond \delta_1$
which is the required $2$-morphism $\nu$.  
\subsection We have then to define  for $j=0,\ldots,n$ an
isomorphism $\beta _j:\ s_j^*\L_{n+1}\isom\L_{n}$. The line bundle $s_j^*\L$ is defined by the
morphism $l\rond p_{n+1}\rond s_j$ and $p_{n+1}\rond s_j$ is associated
to the canonical inclusion $\Delta_0\hookrightarrow\D_n$ which means $p_{n+1}\rond s_j=p_n$.
Therefore, $s_j^*\L_{n+1}=\L_n$ and we define $\beta_j$ to be the identity. 
\subsection We have to show that the data $\L_\bullet,\alpha _j,\beta_j,\ j\geq 0$ defines a line
bundle on the simplicial stack $\D_\bullet$ as explained in (\ref{rel-iso}). Notice that the fact
that the definition of the $\beta_j$'s is compatible with the relations \ref{B} is tautological
($\beta_j$ is the identity on the relevant $\L_n$). 

\subsection Relation \ref{A}: in terms of $l$,
this relation means the following. We have the 2 stricltly commutative diagrams  $$\alpha
_i\rond(\delta_i*\alpha _j):\ l\rond p_n\rond \delta_j\rond \delta_i\FFhd{\delta_i*\alpha
_j}{}l\rond p_{n+1}\rond\delta_i\FFhd{\alpha _i}{}l\rond p_{n+2}$$ diagrams $$\matrix{
\D_{n+2}&\fhd{\delta_i}{}&\D_{n+1}&\fhd{\delta_j}{}&\D_n\cr &\llap{$\scp
p_{n+2}$}\searrow&\fvb{}{\scp p_{n+1}}&\swarrow\rlap{$\scp p_n$}{}\cr &&\D&\fhd{l}{}BG_m\cr }\
{\rm and}\ \matrix{ \D_{n+2}&\fhd{\delta_{j+1}}{}&\D_{n+1}&\fhd{\delta_i}{}&\D_n\cr
&\llap{$\scp p_{n+2}$}\searrow&\fvb{}{\scp p_{n+1}}&\swarrow\rlap{$\scp p_n$}{}\cr
&&\D&\fhd{l}{}BG_m\cr
}$$ 
inducing the two
$2$-morphisms
$$ \alpha _i\rond(\alpha _j*\delta_i):\ l\rond p_n\rond \delta_j\rond
\delta_i\FFhd{\alpha _j*\delta_i}{}l\rond p_{n+1}\rond\delta_i\FFhd{\alpha _i}{}l\rond p_{n+2}$$
and
$$ \alpha _{j+1}\rond(\alpha _i*\delta_{j+1}):\ l\rond p_n\rond \delta_i\rond
\delta_{j+1}\FFhd{\alpha _i*\delta_{j+1}}{}l\rond p_{n+1}\rond\delta_{j+1}\FFhd{\alpha
_{j+1}}{}l\rond p_{n+2}.$$  The relation \ref{A} means exactly the equality 
$$\alpha _i\rond(\alpha _j*\delta_i)=\alpha _{j+1}\rond(\alpha _i*\delta_{j+1}),
\ i\leq j.\leqno(\ref{A}')$$

If $j=0$, the relation \ref{A}' is just by definition of $\alpha _j$ the condition \ref{lineA}
(see remark \ref{rem-line}).

If $j>0$,  both the 2 isomorphisms $\alpha_j$ and $\alpha_{j+1}$ are the relevant identity and
the relation \ref{A}' is tautological.

\subsection Relation \ref{C}: the only non tautological relation in (\ref{C}) is the associated
to the equality $s_0\rond \delta_0=1$ in $(\Delta)$ which means as before that $\alpha
_0*\delta_0$ is the identity functor of $l\rond p_n=l\rond p_n\rond\delta_0\rond s_0$. But, this
is exactly the meaning of the relation \ref{lineB} (see remark \ref{rem-line}).

\section{Application to the Picard groups of some moduli spaces}  Let us chose 3 integers
$r,s,d$ such that
$$r\geq 2\  {\rm and}\ s\mid r\mid ds.$$ 
If $G$ is the group $\SL_r/\mu_s$ we denote as in [B-L-S] by
$\M_G(d)$ the moduli stack of $G$-bundles on $X$ and by $\M_{\SL_r}(d)$ the moduli stack of rank
$r$ vector bundles and determinant $\O(d.x)$. If $r=s$ (i.e. $G={\bf PSL}_r$), the natural morphism
of algebraic stacks $$\pi:\ \M_{\SL_r}(d)\ra\M_G(d)$$ is a $\J_r$-torsor (see the corollary of
proposition 2 of [Gr] for instance). Let me  explain how to deal with the general case.
\subsection Let $E$ be a a  rank
$r$ vector bundle on $X_S$ endowed with an isomorphism $\tau;\ D^{r/s}\isom\det(E)$ 
where $D$ is some line bundle. Let me define the $\SL_r/\mu_s$-bundle $\pi(E)$ associated to  $E$ 
(more precisely to the pair $(E,\tau)$).
\th Definition
\enonce An $s$-trivialization of $E$ on the
\'etale neighborhood $T\ra X_S$ is a triple $(M,\alpha ,\sigma)$ where $\alpha:\ D\isom M^s$ is 
an isomorphism ($M$ is a line bundle on $T$); $\sigma:\ M^{\oplus r}\isom E_T$ is an isomorphism;
$\det(\sigma)\rond \alpha^{r/s}:\ D^{r/s}\isom\det(E)$ is equal to $\tau$. 

Two  $s$-trivializations $(M,\alpha ,\sigma)$ and $(M',\alpha' ,\sigma')$ of $E$ will be said
equivalent if there exists an isomorphism $\iota: M\isom M'$ such that $\iota^s\rond\alpha
=\alpha'$.
\endth

 The principal homogeneous space 
$$T\longmapsto \{\hbox{equivalence classes of $s$-trivializations of } E_T\}$$ 
  defines the
$\SL_r/\mu_s$-bundle $\pi(E)$\footnote{${}^\dagger$}{\sevenrm We see here a $\scp G$-bundle as a
formal homogeneous space under $\scp G$.}. Now, the construction is obviously functorial and
therefore defines the morphsim $\pi:\ \M_{\SL_r}(d)\ra \M_G(d)$ (observe that an object $E$ of
$\M_{\SL_r}(d)$ has determinant $\O({ds\over r}.x)^{r/s}$). Let $L$ be a line bundle and
$(M,\alpha,\tau)$ an $s$-trivialization of $E_T$. Then, $(M\otimes L,\alpha\otimes
\Id_{L^s},\sigma\otimes \Id_L)$ is an $s$-trivialization of $E\otimes L$ (which has determinant
($D\otimes L^s)^{r/s}$). This shows that there exists a canonical functorial isomorphism
$$\pi(E)\isom\pi(E\otimes L)\leqno{(\formule)}.$$\label{iso-pro} In particular, 
$\pi$ is 
$\J_s$-equivariant.

\th Lemma \enonce  The natural morphism of algebraic stacks $\pi:\ \M_{\SL_r}(d)\ra\M_G(d)$
is a $\J_s$-torsor. \endth
\label{pi-tors} {\it Proof}: let $E,E'$ be two rank $r$ vector
bundles on $X_S$ (with determinant equal to $\O(d.x)$) and let $\iota:\ \p(E)\isom\p(E)'$ an
isomorphism. As in the proof of the lemma 13.4 of [B-L-S], we have the exact sequence of sets
$$1\ra\mu_s\ra {\rm Isom}(E,E')\ra{\rm Isom}(\p(E),\p(E)')\fhd{\pi_{E,E'}}{} H^1_{\rm
\acute et}(X_S,\mu_s).$$ Let $L$ be a $\mu_s$-torsor such that $\pi_{E,E'}(\iota)=[L]$. Then,
$\p(E\otimes L)$ is canonically equal to $\p(E)$ and $\pi_{E\otimes L,E'}=0$ and $\iota$ is
induced by an isomorphism $E\otimes L\isom E'$ well defined up to multiplication by $\mu_s$. The
lemma follows.\cqfd 

\subsection Let $\U$ be the universal bundle
on  $X\times \M_{\SL_r}(d)$. We would like to know which power of the determinant bundle $\d=(\det
R\Gamma\U)^{-1}$
 on $\M_{\SL_r}(d)$ descends to $\M_G(d)$. As in I.3 of [B-L-S], the rank $r$ bundle
$\F=\L^{\oplus(r-1)}\oplus\L^{1-r}(d.x)$ on $X\times\J(X)$ has determinant $\O(d.x)$ and therefore
defines a morphism $$f:\ \J(X)=JX\times BG_m\ra\M_{\SL_r}(d)$$ which is $\J_s$-equivariant.

The
vector bundle $\F'=\O^{\oplus (r-1)}\oplus\L^{-r/s}(d.x)$ on $X\times\J(X)$ has determinant
$[\L^{-1}({ds\over r}.x)]^{r/s}$. The
$G$-bundle 
$\pi(\F')$ on $X\times\J(X)$ defines a morphism $f':\ \J\ra\M_G(d)$. The relation
$\L\otimes (\Id_X\times s_\J)^**(\F')=\F$ and (\ref{iso-pro}) gives an isomorphism
$\pi(\F)=(\Id_X\times s_\J)^*\pi(\F')$ which means that the diagram
 
$$\matrix{
\J(X)&\fhd{f}{}&\M_{\SL_r}(d)\cr
\fvb{}{s_\J}&&\fvb{}{\pi}\cr
\J(X)&\fhd{f'}{}&\M_G(d)\cr}\leqno{(\formule)}$$\label{dia-rest} is $2$-commutative.
Exactly as in I.3 of [B-L-S], let me prove the
\th Lemma
\enonce The line bundle $f^*\d^k$ on $\J(X)$ descends through $s_\J$ if and only if $k$
multiple of $s/(s,r/s)$.
\endth\label{lem-rest}
{\it Proof}: let
$\chi=r(g-1)-d$ be the opposite of the Euler characteristic of ($\k$-)points of $\M_{\SL_r}(d)$. By
(\ref{iso-det}), one has an isomorphism
$f^*\d^k\isom\Theta^{kr(r-1)}\otimes\T^{k\chi}.$
The theory of Mumford groups says that $\Theta^{kr(r-1)}$ descends through $s_J$ if and only if $k$
is a multiple of $s/(s,r/s)$. The line bundle $\T^{k\chi}$ on $BG_m$ descends through $s_{BG_m}$ if
and only if $k\chi$ is a multiple of $s$. The lemma follows from the above isomorphism and from the
observation that the condition
$s\mid r\mid ds$ forces $s\chi$ to be a multiple of $s$.\cqfd

\rem\label{rema} If $g=0$, the jacobian $J$ is a point and the condition on
$\Theta$ is empty. The only condition is in this case being $k\chi$ multiple
of
$s$.

Let me recall that $\d$ is the determinant bundle on $\M_{\SL_r}(d)$ and $G=\SL_r/\mu_s$. 
\th Theorem
\enonce Assume that the characteristic of $\k$ is $0$. The integers $k$ such
that $\d^k$ descends to $\M_G(d)$ are the multiple of $s/(s,r/s)$.
\endth\label{theo-pic}

{Proof}: by lemma \ref{lem-rest} and diagram (\ref{dia-rest}), we just have to proving that $\d^k$
efectively descends where $k=s/(s,r/s)$. By theorem \ref{theo} and lemma \ref{pi-tors}, this means
exactly that $\d^k$ has a $\J_s$-linearization. We know by lemma \ref{lem-rest} that the pull-back
$f^*\d^k$ has such a linearization.
\th Lemma
\enonce The pull-back morphism $\Pic(\J_s\times\M_{SL_r}(d))\ra\Pic(\J_s\times\J(X))$ is injective.
\endth
{\it Proof}: by lemma \ref{torseur}, one is reduced to prove that the natural morphism
$$\Pic(B\mu_s\times\M_{SL_r}(d))\ra\Pic(B\mu_s\times\J(X))$$ is injective. Let  $\D$ be any stack.
The canonical morphism $\D\ra\D\times B\mu_s$ is a $\mu_s$-torsor (with the trivial action of
$\mu_s$ on $\D$). By theorem \ref{theo}, one has the equality $$\Pic(\D\times
B\mu_s)=\Pic^{\mu_s}(\D).$$ Assume further that 
 $H^0(\D,\O)=\k$. The later group  is then canonically isomorphic to
$$\Pic(\D)\times\Hom(\mu_s,G_m)=\Pic(\D)\times\Pic(B\mu_s).$$ All in all, we get a functorial
isomorphism $$\iota_\D:\ \Pic(\D\times
B\mu_s)\isom\Pic(\D)\times\Pic(B\mu_s).\leqno(\formule)\label{pic-pro}$$ By [L-S], the Picard group
of $\M_{\SL_r}(d)$ is the free abelian group $\Z.\d$ and the formula (\ref{iso-det})
proves that  $$f^*:\ \Pic(\M_{\SL_r}(d))\ra\Pic(\J(X))$$ is an injection. The lemma folows from the
commutative diagram
$$\diagram{
\Pic(\M_{\SL_r}(d))\times\Pic(B\mu_s)&\hookrightarrow&\Pic(\J(X))\times\Pic(B\mu_s)\cr
\fvb{\iota_\M}{\wr}&&\fvb{\iota_\J}{\wr}\cr
\Pic(\M_{\SL_r}(d)\times B\mu_s)&\ra&\Pic(\J(X)\times B\mu_s)\cr
}$$\cqfd 
 Let $\H$ (resp. $\H_\J$) be the line bundle on $\J_s\times\M_{\SL_r}(d)$ (resp.
$\J_s\times\J(X)$  
$$\H=\hom(m_\M^*\d^k,pr_2^*\d^k)\ {\rm resp.}\
\H_\J=\hom(m_\M^*f^*\d^k,pr_2^*f^*\d^k).$$
Let us chose a $\J_s$-linearization $\lambda_\J$ of $f^*\d^k$. It defines a trivialization of the line
bundle $\H_\J$.
 The equivariance of $f$ implies (cf. \ref{equi})   that
there exists a (compatible) $2$-morphism 
$$q:\ m_\M\rond(\Id\times f)\Rightarrow f\rond m_\J$$
making the diagram 
$$ \matrix{ \J_s\times\J(X)&\fhd{m_\J}{}&\J(X)\cr \fvb{{\rm Id}\times
f}{}&&\fvb{f}{}\cr \J_s\times\M_{SL_r}(d)&\fhd{m_\M}{}&\M_{SL_r}(d)\cr }$$  
$2$-commutative. The $2$-morphism $q$ defines an isomorphism from the pull-back  $m_\M^*\d^k$ on
$\J_s\times\J(X)$ to $m_\J^*(f^*\d^k)$. The pull-back of $pr_2^*\d^k$ on $\J_s\times\J(X)$ is
tautologically isomorphic to $pr_2^*(f^*\d^k)$. The
preceding isomorphisms induce an isomorphism $$(\Id\times f)^*\H\isom\H_J.$$  The later line
bundle  being trivial, so is $(\Id\times f)^*\H$. The lemma above proves therefore that $\H$
istself is {\it trivial}  . Each ($\k$-)point $j$ of $J_s$ defines a morphism
$\M_{SL_r}(d)\ra\J_s\times\M_{SL_r}(d)$ (resp. $\J(X)\ra\J_s\times\J(X)$); let me denote by $\H_j$
(resp. $f^*\H_j$) the pull-back of $\H$ (resp. $(\Id\times f)^*\H)$ by this morphism. The pull-back
morphism $$H^0(\J_s\times\M_{\SL_r}(d),\H)\ra H^0(\J_s\times\J(X),(\Id\times f)^*\H)$$ can be
identified to the direct sum $$\oplus_{j\in\J_s(\k)}H^0(\M_{\SL_r}(d),\H_j)\ra
H^0(\J(X),f^*\H_j).$$ 
Because 
$$H^0(\M_{\SL_r}(d),\O)=H^0(\J(X),\O)=\k\leqno{(\formule)},\label{cons}$$ this morphism is a direct
sum of non-zero morphisms of vector spaces of dimension $1$ and therefore an isomorphism. In
particular, a linearization $\lambda_\J$ of $f^*\d^k$ defines canonicaly an isomophism $$\lambda_\M:\
m_\M^*\d^k\isom pr_2^*\d^k$$ such that $(\Id\times f)^*\lambda_\M=\lambda_\J$. 

Explicitely, $\lambda_\M$ is
characterized as follows: let $x$ be an object of $\M_{\SL_r}(d)$ over a connected scheme $S$ and
$g$ an object of $\J_s(S)=\J_s(\k)$. The preceding dicussion means that the functorial
isomorphisms $$\lambda_\M(g,x):\ \d^k_{g.x}\isom\d^k_x$$ are determined when $x$ lies in the essential
image of $f$. In this case, let us chose an isomorphism $f(x')\isom x$ (inducing an isomorphism 
$g.f(x')\isom g.x$).
Then, the diagram of isomorphisms of line bundles on $S$
$$\matrix{
L'_{x'}=L_{f(x')}&\ra&L_x\cr
\fvh{\lambda_\J(g,x')}{}&&&\llap{$\scp\lambda_\M(g,x)$}\nwarrow&\cr
L'_{g.x'}=L_{f(g.x')}&\fhd{q_{g,x'}}{}&L_{g.f(x')}&\ra&L_{g.x}\cr
}$$ is commutative (where $L=\d^k$ and $L'=f^*\d^k$).

 Now, the pull-back of $\lambda_\M$ on
$\J_s\times\J(X)$ satisfies  conditions \ref{lineA} and \ref{lineB}. Using (\ref{cons}) and the
equivariance of of $f$ as above, this shows that $\lambda_\M$ is a linearization. For instance,
keeping the notation above, let us check the condition \ref{lineB}. We have to check that the
isomorphism $\iota$ of $L$ induced by $\epsilon$ is the identity. As above, it is enough to check
that on $\J(X)$.  By definition, with a slight abuse of notations,  the diagrams
$$\matrix{
L'_{x'}=L_{f(x')}&\ra&L_x\cr
\fvh{\lambda_\J(1,x')}{}&&&\llap{$\scp\lambda_\M(1,x)$}\nwarrow&\cr
L'_{1.x'}=L_{f(1.x')}&\fhd{q_{1,x'}}{}&L_{1.f(x')}&\ra&L_{1.x}\cr
}\kern .5cm {\rm and}\kern .5cm
\matrix{
L_x&\fhd{\iota}{}&L_x\cr
&\llap{$\scp\lambda_\M(1,x)$}\nwarrow&\fvh{\epsilon(x)}{}\cr
&&L_{1.x}\cr
}$$ 
Because $\lambda_\J$ is a linearization, condition \ref{lineB} gives the commutative diagram
$$\matrix{
L'_{x'}&=&L'_{x'}\cr
&\llap{$\scp\lambda_\J(1,x')$}\nwarrow&\fvh{}{\epsilon'(x')}\cr
&&L_{1.x'}\cr
}$$ showing that the equality $\iota=\Id$ remains to prove the commutativity of the diagram
$$\matrix{
L_{f(1.x')}&\fhd{\epsilon'}{}&L_{f(x')}\cr
\fvb{q_{1,x'}}{}&&\parallel\cr
L_{1.f(x')}&\fhd{\epsilon}{}&L_{f(x')}\cr
}
$$ 
But this follows from the commutativity of the diagram
$$\matrix{
f(1.x')&\fhd{\epsilon'}{}&f(x')\cr
\fvb{q_{1,x'}}{}&&\parallel\cr
1.f(x')&\fhd{\epsilon}{}&f(x')\cr
}
$$ which is by definition the fact that $q$ is compatible to $\epsilon$ as required in
(\ref{equi}). One would check condition \ref{lineA} in an analogous way.\cqfd 
\rem\label{rem}In the case $g=0$, the condition is an in remark \ref{rema}.

\rem This linearization can be certainly also deduced from a careful
analysis of the first section of [Fa], but the method above seems simpler.
\bigskip

\centerline{\bf References}

\medskip
[B-L-S] A.Beauville, Y. Laszlo, C. Sorger, {\it The Picard group of the moduli of $G$-bundles on
a curve}, preprint alg-geom/9608002. 
\medskip
[Br] L. Breen, {\it Bitorseurs et cohomologie non ab\'elienne
}, in The Grothendieck Festschrift I, Progr. Math. {\bf 86},
Birkh\"auser (1990), 401-476.
\medskip
[De] P. Deligne, {\it Th\'eorie de Hodhe III}, Publ. Math.
I.H.E.S. {\bf 44} (1974), 5-78.
\medskip
[Fa] G. Faltings, {\it Stable $G$-bundles and Projective Connections}, JAG {\bf
2} (1993), 507-568.
\medskip
[Gr] A. Grothendieck, {\it G\'eom\'etrie formelle et g\'eom\'etrie
alg\'ebrique}, FGA, S\'em. Bourbaki {\bf 182} (1958/59), 1-25.
\medskip
[L-M] G. Laumon, L. Moret-Bailly, {\it Champs alg\'ebriques}, preprint
Universit\'e Paris-Sud (1992).
\medskip
[L-S] Y. Laszlo, C. Sorger, {\it The line bundles on the moduli of parabolic $G$-bundles over
curves and their sections}, preprint alg-geom/9507002.  
\bigskip
[McL] S. Mac Lane, {\it Categories for the working mathematician}, GTM {\bf  5}, Springer-Verlag (1971).
\bigskip

\hfill\hbox to 5cm{\hfill Y. {\pc LASZLO}\hfill}
\smallskip

\hfill\hbox to 5cm{\hfill DMI - \'Ecole Normale Sup\'erieure\hfill}
\smallskip

\hfill\hbox to 5cm{\hfill ( URA 762 du CNRS )\hfill}
\smallskip

\hfill\hbox to 5cm{\hfill 45 rue d'Ulm\hfill}
\smallskip

\hfill\hbox to 5cm{\hfill F-75230 {\pc PARIS} Cedex 05\hfill}
\smallskip

\end